\documentclass[a4paper,english,twoside,12pt]{smfart}

\usepackage{smfenum}
\usepackage{smfthm}
\usepackage{amssymb}
\usepackage{mathrsfs}
\input{xypic}

\def\CC{{\mathbb C}}

\def\NN{{\mathbb N}}

\def\RR{{\mathbb R}}
\def\ZZ{{\mathbb Z}}

\def\D{{\rm D}}
\def\E{{\rm E}}
\def\F{{\rm F}}
\def\G{{\rm G}}

\def\I{{\rm I}}
\def\J{{\rm J}}
\def\K{{\rm K}}
\def\L{{\rm L}}
\def\M{{\rm M}}
\def\N{{\rm N}}
\def\P{{\rm P}}
\def\Q{{\rm Q}}
\def\R{{\rm R}}

\def\U{{\rm U}}
\def\V{{\rm V}}
\def\W{{\rm W}}

\def\Cc{\mathscr{C}}
\def\Dd{\mathscr{D}}

\def\Hh{\mathscr{H}}
\def\Ii{\mathscr{I}}

\def\Nn{\mathscr{N}}
\def\Oo{\mathscr{O}}

\def\Rr{\mathscr{R}}
\def\Ss{\mathscr{S}}
\def\Tt{\mathscr{T}}
\def\Vv{\mathscr{V}}

\def\d{\delta}

\def\f{\rightarrow}

\def\h{\varphi}

\def\k{\kappa}
\def\l{\lambda}

\def\p{\mathfrak{p}}

\def\s{\sigma}

\def\lp{\langle}
\def\rp{\rangle}
\def\>{\geqslant}
\def\<{\leqslant}

\def\Hom{{\rm Hom}}
\def\End{{\rm End}}

\def\GL{{\rm GL}}

\def\Ind{{\rm Ind}}

\def\Mat{{\rm M}}

\def\mult#1{{#1}^{\times}}

\def\ss{\Ss}
\def\supp{{\rm supp}}
\def\Irr{{\rm Irr}}

\def\SS{{\rm S}}

\def\k{k}

\author{Vincent S\'echerre}
\address{Institut de Math\'ematiques de Luminy\\
CNRS UMR $6206$\\
Universit\'e de Luminy\\
163 avenue de Luminy\\
$13288$ Marseille Cedex $9$\\
France}
\email{secherre@iml.univ-mrs.fr}

\numberwithin{equation}{section} 
\theoremstyle{plain}

\title[Proof of the Tadi\'c conjecture U0]
{Proof of the Tadi\'c conjecture U0 on the unitary dual of $\GL_m(\D)$}

\begin{abstract}
Let $\F$ be a non-Archimedean local field of characteristic $0$,
and let $\D$ be a finite dimensional central division algebra over
$\F$.
We prove that any unitary irreducible representation of a Levi 
subgroup of $\GL_m(\D)$, with $m\>1$, induces irreducibly to
$\GL_m(\D)$.
This ends the classification of the unit\-ary dual of $\GL_m(\D)$
initiated by Tadi\'c.
\end{abstract}

\thanks
{The research for this paper was partially supported by EPSRC grant
GR/T21714/01}

\begin{document}

\maketitle

\section*{Introduction}

Let $\F$ be a non-Archimedean locally compact non-discrete field of
characteristic zero (that is, a finite extension of the field of
$p$-adic numbers for some prime number $p$) and let $\D$ be a
finite-dimensional central division algebra over $\F$.
In \cite{Tadic2}, Tadi\'c gave a conjectural class\-if\-ication
of the unitary dual of $\GL_m(\D)$ based on five statements
U0,\ldots,U4; in the same article, he proved U3 and U4.
In \cite{BaRe}, Badulescu and Renard proved U1, and it is known that
U0 and U1 together imply U2. 
In this paper, we prove the remaining conjecture U0, which asserts
that any unit\-ary irr\-educible representation of a Levi subgroup of 
$\GL_m(\D)$, with $m\>1$, induces irreducibly to $\GL_m(\D)$.
The proof is based on Bushnell-Kutzko's the\-ory of types (see
\cite{BK1}), and more precisely on their theory of covers, which
allows one to com\-pare parabolic induction in $\GL_m(\D)$ with 
parabolic induction in affine Hecke algebras. 

The proof consists of reducing to the case where $\D$ is commutative,
for which the result is already known (Bernstein \cite{BernsteinPinv}, 
see Theorem \ref{bouloulou} below).
This can be done by using particular types of $\GL_m(\D)$, the
so-called Bushnell-Kutzko simple types (see \cite{BK,VS3}).
Their Hecke algebras are well known and isomorphic to affine
Hecke-Iwahori algebras, which allows one to transport our induction
problem, {\it via} the Hecke algebra isomorphisms of \cite{VS3}, to a
very special case, in which the conjecture is known to be true. 
This method has been already used in \cite{BHK,BK}. 

The proof can be decomposed into three parts.
In the first part (\S\ref{OnDenoting}), we reduce to the case where
the unitary irreducible representation of the Levi subgroup is simple
in the sense of \cite{BK1}: the elements of its cuspidal support are
unramified twists of a single cuspidal irreducible representation of
$\GL_{k}(\D)$, where $k$ is a divisor of $m$.
This special case of the conjecture is denoted by S0.
In the second part, we translate the problem in
terms of induction of mod\-ules over Hecke algebras.
More precisely, we reduce the proof of S0 to proving 
that, given $r\>1$, any unitary irreducible module over 
the Hecke-Iwahori algebra of a Levi subgroup of $\GL_r(\F)$ induces
irreducibly to the Hecke-Iwahori algebra of $\GL_r(\F)$
(see Proposition \ref{Diderot}).
This step demands the existence of covers for any irreducible
simple representation of $\GL_m(\D)$.
Such covers have been constructed in \cite{VS1,VS2,VS3,VS4}.
The last part of the proof consists of proving Proposition
\ref{Diderot} (see above).
This step is based on a result of Barbasch-Moy (see \cite{BM1,BM2}) 
which asserts that the functor of Iwahori-invariant vectors induces a
one-to-one correspondence between: 
\begin{itemize}
\item[(1)]
unitary irreducible representations of $\GL_r(\F)$ having a non-zero
vector invariant under an Iwahori subgroup;
\item[(2)]
unitary irreducible modules over the Hecke-Iwahori algebra of
$\GL_r(\F)$.
\end{itemize}

In the last section of this article, we determine the unramified
characters $\chi$ of $\GL_m(\D)$ for which the parabolically induced
representation $\Pi(\chi)=\rho\times\rho\chi$, where $\rho$ is a fixed
cuspidal irreducible representation of $\GL_m(\D)$, is reducible.
Unlike \cite{Tadic2}, our result does not refer to the
Jacquet-Langlands correspondence.
This answers a question of J.\ Bernstein and A.\ M\'\i nguez.
Here again, we reduce to the case where $\D$ is commutative, for which
the reducibility points are known to be $\chi=|\det|_\F$ and 
$\chi=|\det|_\F^{-1}$, where $|\ |_\F$ denotes the normalized absolute
value of $\F$.
However, in the division algebra case, the reducibility points $\chi$
depend on the cuspidal representation $\rho$ (see Theorem
\ref{MainTheo4}).

\Subsection*{Acknowledgments}

I would like to thank I.\ Badulescu, D.\ Barbasch, J.\ Bernstein, 
A.\ Minguez, G.\ Mui\'c and S.\ Stevens for helpful conversations. 
Special thanks to Ernst-Wilhelm Zink for stimulating discussions
during my stay at Humboldt University, to Shaun Stevens for linguistic
comments and to University of East Anglia for its hospitality.
I am also grateful to Guy Henniart for suggesting that I look at
Tadi\'c's conjecture, and for his comments on earlier drafts of the
manuscript.


\section{Notations and preliminaries}

In this section, we fix some notations and recall some well-known facts.
The reader may refer to \cite{Tadic2} for more details.

\subsection{}

Let $\F$ be a non-Archimedean locally compact non-discrete field of
characteristic $0$, and let $\D$ be a finite-dimensional central
division algebra over $\F$.
For any integer $m\>1$, we denote by $\M_m(\D)$ the $\F$-algebra 
of $m\times m$ matrices with coefficients in $\D$ and by 
$\G_m=\GL_m(\D)$ the group of its invertible elements.
For convenience, $\G_0$ will denote the trivial group.

Let $\N_m$ be the reduced norm of $\M_m(\D)$ over $\F$ and let 
$|\ |_{\F}$ be the normalized absolute value of $\F$.
The map $g\mapsto|\N_m(g)|_\F$ is a continuous group homomorphism
from $\G_m$ to the multiplicative group $\mult\CC$ of the field of
complex numbers, which we simply denote by $\nu$.

If $\rho$ is a representation and $\chi$ a character of $\G_m$ for
some $m$, we denote by $\rho\chi$ (or equivalently by $\chi\rho$) the
twisted representation $g\mapsto\chi(g)\rho(g)$.

We denote by $\NN$ the set of non-negative integers.
If $\SS$ is a set, a {\it multiset} on $\SS$ is a finitely supported
function from $\SS$ to $\NN$.
It can be thought as an unordered finite family of elements of $\SS$.
For $n\>0$ and $x_{i}\in\SS$ with $1\<i\<n$, we denote by:
\begin{equation*}
\label{MultisetNotation}
(x_{1},\ldots,x_{n})
\end{equation*}
the multiset whose value on $x\in\SS$ is the number of integers
$1\<i\<n$ such that $x_{i}=x$.
The integer $n$ is called the {\it size} of this multiset.
We denote by $\M(\SS)$ the set of all multisets on $\SS$.
It is naturally endowed with a structure of commutative semigroup.

\subsection{}
\label{UZero}

For $m\>0$, we denote by $\Irr_m$ the set of all classes of
irreducible representations of $\G_m$, by $\Rr_m$ the category of
smooth complex representations of finite length of $\G_m$ and by
$\R_m$ the Grothendieck group of $\Rr_m$, which is a free $\ZZ$-module
with basis $\Irr_m$.
In particular, $\Irr_0$ is reduced to a single element and $\R_0$ is
isomorphic to $\ZZ$.
For $\s\in\Irr_{m}$, we set $\deg(\s)=m$, which we call the 
{\it degree} of $\s$.
We set:
\begin{equation*}
\R=\bigoplus\limits_{m\>0}\R_m
\end{equation*}
and:
\begin{equation*}
\Irr=\bigcup\limits_{m\>0}\Irr_m.
\end{equation*}
The group $\R$ is a graded free $\ZZ$-module with basis $\Irr$.
Two equivalent irreducible representations will be considered as the
same element of $\Irr$.

Given $m,n\>0$, the (normalized) parabolic induction functor:
\begin{eqnarray*}
\Rr_{m}\times\Rr_{n} & \f & \Rr_{m+n}\\
(\s,\tau) & \mapsto & \s\times\tau
\end{eqnarray*}
induces a map $\R_{m}\times\R_{n}\f\R_{m+n}$.
This map extends to a $\ZZ$-bilinear map $\R\times\R\f\R$, which makes
$\R$ into an associative and commutative graded $\ZZ$-algebra 
(see \cite[\S2.3]{BZ1} and \cite[\S1]{Tadic2}).
The image of $(\s,\tau)\in\R\times\R$ by this map will be still 
denoted by $\s\times\tau$.

We will make no distinction between {\it unitary} and 
{\it unitarizable} irreducible representations, which form a subset of
$\Irr$ denoted by $\Irr^u$ (see \cite[\S2.8]{Cartier}).
Conjecture U0 is the following statement (see \cite[\S6]{Tadic2}):

\medskip

\noindent{\bf (U0)} \ 
{\it Let $\s,\tau\in\Irr^u$ be unitary irreducible representations.  
Then $\s\times\tau\in\Irr$.}

\medskip

Let us recall the following result of Bernstein \cite{BernsteinPinv}.

\begin{theo}[Bernstein]
\label{bouloulou}
Assume that $\D=\F$.
Then {\rm U0} is true.
\end{theo}

\subsection{}
\label{multiSet}

Let $\Cc$ be the set of all cuspidal representations in $\Irr$.
Given $\s\in\Irr$, there exists a unique multiset:
\begin{equation*}
(\rho_1,\ldots,\rho_n)\in\M(\Cc)
\end{equation*}
such that $\s$ is a subquotient of the induced representation
$\rho_1\times\ldots\times\rho_n$ (see \cite[\S2]{BZ1}).
This multiset is denoted by $\supp(\s)$, and it is called 
the (cuspidal) {\it support} of $\s$.
This defines a surjective map $\supp:\Irr\f\M(\Cc)$, which extends to
$\R$ by linearity.

Given $\s,\tau\in\Irr$, we have:
\begin{equation}
\label{MultSupp}
\supp(\s\times\tau)=\supp(\s)+\supp(\tau).
\end{equation}

More generally, let $\M$ be a Levi subgroup of $\G_m$. 
This subgroup is equal, up to conjugacy by an element of $\G_m$, 
to $\G_{m_1}\times\ldots\times\G_{m_l}$, where the
$m_i$ are positive integers such that $m_1+\ldots+m_l=m$.
Any irreducible representation of $\M$ is of the form 
$\s_1\otimes\ldots\otimes\s_{l}$ with $\s_i\in\Irr_{m_i}$.
The (cuspidal) {\it support} of such a representation is the
sum of the $\supp(\s_i)$, for $1\<i\<l$.

\subsection{}
\label{SerieDiscrete}

Let $\rho\in\Cc$ be a cuspidal irreducible representation, and let $m$
denote the degree of $\rho$.
Let $d$ be the reduced degree of $\D$ over $\F$, that is, 
the square root of the dimension of $\D$ over $\F$.
By the Jacquet-Langlands correspondence (see \cite{DKV}) 
one associates to $\rho$ an essentially square int\-egrable
representation $\s$ of the group $\GL_{md}(\F)$.
The classification of the discrete series of $\GL_{md}(\F)$ 
(see \cite{Ze}) gives us a unique positive integer $b$ dividing 
$md$ and a unique cuspidal irreducible representation $\tau$ of 
$\GL_{md/b}(\F)$ such that $\s$ is a quotient of the induced
representation
$\tau\times\mu\tau\times\ldots\times\mu^{b-1}\tau$, where
$\mu:g\mapsto|\det(g)|_\F$ denotes the analogue of $\nu$ for the 
group $\GL_{md/b}(\F)$.
We denote this integer by $b(\rho)$, and we set:
\begin{equation*}
\nu_{\rho}=\nu^{b(\rho)}.
\end{equation*}

Let $\Dd$ be the set of all essentially square integrable 
representations in $\Irr$.
It is parametrized by means of cuspidal irreducible representations  
as follows.
For any $\rho\in\Cc$ and any positive integer $n$, the induced
representation: 
\begin{equation*}
\label{segmentnormalise}
\nu_{\rho}^{(n-1)/2}\rho\times\nu_{\rho}^{-1+(n-1)/2}\rho
\times\ldots\times\nu_{\rho}^{-(n-1)/2}\rho
\end{equation*}
has a unique essentially square integrable quotient,
which we denote by $\d(\rho,n)$.
The map $\Cc\times\NN^*\f\Dd$ obtained this way is a bijection 
(see \cite{Tadic2,Ze}), where $\NN^*$ denotes the set of positive
integers.

Let $\Cc^{u}$ (resp. $\Dd^{u}$) be the set of all unitary
representations in $\Cc$ (resp. in $\Dd$).
Then $\d(\rho,n)$ is unitary if and only if $\rho$ is.
In other words, the image of $\Cc^{u}\times\NN^*$ by the map above is
$\Dd^{u}$.

\subsection{}
\label{JJR}

Let $\Tt$ be the set of all essentially tempered representations 
in $\Irr$ and let $\Tt^{u}$ be the set of all tempered representations
in $\Tt$. 
Given $\tau\in\Tt$, there exists a unique real number $e(\tau)\in\RR$, 
which we call the {\it exponent} of $\tau$, 
such that $\nu^{-e(\tau)}\tau$ is tempered.
The map:
\begin{equation}
\label{Tohoku}
(\d_1,\ldots,\d_k)\mapsto\d_1\times\ldots\times\d_k
\end{equation}
induces a bijective correspondence from $\M(\Dd^{u})$ 
onto $\Tt^{u}$ (see \cite[B.2.$d$]{DKV}). 

Given $d=(\d_1,\ldots,\d_k)\in\M(\Dd)$, the fibers of the map
$i\mapsto e(\d_i)$ decompose $\{1,2,\ldots,k\}$
into a finite disjoint union $\I_1\cup\ldots\cup\I_{l}$.
For $1\<i\<l$, we denote by $\tau_i$ the product of the $\d_j$ for
$j\in\I_i$.
Each $\tau_i$ is essentially tempered.
Let us choose an ordering such that:
\begin{equation*}
\label{stdorder}
e(\tau_1)\>\ldots\>e(\tau_l).
\end{equation*}
Then the induced representation $\tau_1\times\ldots\times\tau_l$
has a unique irreducible quotient, which we denote by $\L(d)$.
This representation depends only on $d$ and not on the ordering 
of the $\tau_i$, and the map $d\mapsto\L(d)$ is a bijection from
$\M(\Dd)$ to $\Irr$.

\subsection{}
\label{hermite}

Given $\s\in\Irr$, we denote by $\s^{\vee}$ the contragredient 
representation of $\s$ and by $\overline{\s}$ its complex conjugate 
representation, that is, the representation obtained 
by making $\CC$ act on the space of $\s$ by
$(\l,v)\mapsto\overline{\l}v$.
The representation:
\begin{equation*}
\s^+=\overline{\s^{\vee}}
\end{equation*}
is called the {\it Hermitian contragredient} of $\s$, and $\s$ is
said to be {\it Hermitian} if it is equivalent to its Hermitian
contragredient.
Since this is equivalent to the existence of a non-degenerate 
invariant Hermitian form on the space of $\s$, 
any unitary irreducible representation is Hermitian.

Given $d\in\M(\Dd)$, we denote by $d^+$ the multiset on 
$\Dd$ whose elements are the Hermitian contragredients 
of the elements of $d$.
Then (see \cite[\S2]{Tadic2}) we have:
\begin{equation*}
\L(d)^{+}=\L(d^+).
\end{equation*}
Thus $\L(d)$ is Hermitian if and only if $d^{+}=d$.
Note that, for $\d\in\Dd$, the exponent of $\d^+$ is $-e(\d)$. 

\begin{lemm}
\label{MuicSaidIt}
Let $\s,\tau\in\Irr$ be Hermitian representations such that 
$\s\times\tau$ is irreducible and unitary.
Then $\s$ and $\tau$ are unitary.
\end{lemm}

\begin{proof}
This is a standard result.
The Hermitian forms on the spaces of $\s$ and $\tau$ induce 
a Hermitian form $h$ on the space of $\s\times\tau$.
As $\s\times\tau$ is irreducible, its space can be endowed
with a unique, up to a non-zero real scalar, non-degenerate 
Hermitian form.
Therefore, up to a sign, $h$ is positive definite, and $\s,\tau$ are
unitary (see \cite[\S3(a)]{Tadic3}).
\end{proof}

\subsection{}

Given $\rho\in\Cc$, we set:
\begin{equation*}
\ell(\rho)=\{\nu^{s}\rho\ |\ s\in\CC\}.
\end{equation*}
A {\it line} in $\Cc$ is a subset of $\Cc$ of the form $\ell(\rho)$
for some $\rho\in\Cc$.

\begin{defi}
\begin{enumerate}
\item[(i)]
An irreducible representation $\s$ is said to be {\it simple}
if there exists a line $\ell$ in $\Cc$ such that $\supp(\s)\in\M(\ell)$.
\item[(ii)]
Two representations $\s,\tau\in\Irr$ are said to be {\it aligned}
if $\s\times\tau$ is simple.
\end{enumerate}
\end{defi}

\begin{rema}
Any essentially square integrable irreducible representation is simple.
If two representations $\s,\tau\in\Irr$ are aligned,
then $\s$ and $\tau$ are simple.
In particular, a representation is simple if and only if it is
aligned with itself. 
\end{rema}

The following result is an immediate consequence of \cite{Tadic2}, 
Proposition $2.2$ and Lemma $2.5$.

\begin{prop}
\label{Tad}
Let $d=(\d_1,\ldots,\d_k)$ and $d'=(\d'_1,\ldots,\d'_{k'})$
be in $\M(\Dd)$.
Suppose that, for any $1\<i\<k$ and $1\<j\<k'$, the 
representations $\d_i$ and $\d'_j$ are not aligned.
Then $\L(d)\times\L(d')$ is irreducible and equal to $\L(d+d')$.
\end{prop}

This leads to the following result.

\begin{prop}
\label{Lapounet}
Let $\s\in\Irr$ be an irreducible representation.
\begin{itemize}
\item[(i)]
There is a unique subset $\{\s_1,\ldots,\s_k\}$ of $\Irr$ such
that $\s=\s_1\times\ldots\times\s_k$, and such that $\s_i,\s_j$ are
aligned if and only if $i=j$.
\item[(ii)]
If $\s$ is unitary, then so are the $\s_i$.
\end{itemize}
\end{prop}

\begin{proof}
Let $d\in\M(\Dd)$ be such that $\s=\L(d)$.
The multiset $d$ can be written in a unique way as a sum:
\begin{equation}
\label{SumMD}
d=d_1+\ldots+d_k
\end{equation}
such that two elements of $d$ are aligned if and only if they are
cont\-ain\-ed in the same $d_i$.
Thus, according to Proposition \ref{Tad}, we have:
\begin{equation*}
\L(d)=\L(d_1)\times\ldots\times\L(d_k).
\end{equation*}
The unicity property comes from the unicity of decomposition
(\ref{SumMD}).
Moreover, if $d^+=d$, then $d^+_i=d_i$ for each integer $1\<i\<k$.
Therefore, if $\L(d)$ is Herm\-itian, then so are the $\L(d_i)$.
By Lemma \ref{MuicSaidIt}, if $\L(d)$ is unitary, then 
so are the $\L(d_i)$.
\end{proof}


\section{Theory of types for $\GL_{m}(\D)$}

In order to prove Conjecture U0, we need some material from
Bushnell-Kutzko's theory of types, which we develop in this section.

\subsection{}
\label{RappelHecke}

Let $m$ be a positive integer, and let $\M$ be a Levi subgroup of
$\G=\G_{m}$.
Let $\J$ be a compact open subgroup of $\M$, and let $\tau$ be a
smooth irreducible representation of $\J$ on a complex vector space
$\Vv$.
Let us choose a Haar measure on $\M$ giving measure $1$ to $\J$.
The {\it Hecke alg\-ebra} of $\M$ relative to $(\J,\tau)$,
which we de\-no\-te by $\Hh(\M,\tau)$, is the convolution algebra 
of locally constant and compactly supported func\-tions
$f:\M\f\End_{\CC}(\Vv)$ such that:
\begin{equation*}
f(kgk')=\tau(k)\circ f(g)\circ\tau(k')
\end{equation*}
for any $k,k'\in\J$ and $g\in\M$.
We have a functor:
\begin{equation}
\label{FoncMT}
{\bf M}_{\tau}:\s\mapsto\Hom_{\J}(\tau,\s)
\end{equation}
from the category of smooth complex representations of $\M$ to the
category of right modules over $\Hh(\M,\tau)$. 
It induces a bijection between the classes of irred\-ucible
representations of $\M$ whose restriction to $\J$ contains $\tau$ and
the classes of irred\-ucible right $\Hh(\M,\tau)$-modules.

\subsection{}
\label{StructureUnitaire}

According to \cite[\S4.3]{BK}, the Hecke algebra $\Hh(\M,\tau)$ can
be canonically endowed with an involution $f\mapsto f^*$.
A right module $\V$ over $\Hh(\M,\tau)$ is said to be {\it unitary} 
if there exists a positive definite Hermitian form 
$(x,y)\mapsto\lp x,y\rp$ on $\V$ such that: 
\begin{equation*}
\lp vf,w\rp=\lp v,wf^*\rp
\end{equation*}
for any $v,w\in\V$ and $f\in\Hh(\M,\tau)$.

Note that ${\bf M}_\tau$ preserves unitarity: if an irreducible
representation of $\M$ is unitary, then the irreducible module which
corresponds to it is unitary. 

\subsection{}
\label{RappelTypes}

Let $(\rho_1,\ldots,\rho_k)\in\M(\Cc)$ be a multiset of cuspidal
irreducible representations of $\M$.
The {\it inertial class} of this multiset is the set $\ss$ of all
multisets of the form $(\rho_1\chi_1,\ldots,\rho_k\chi_k)$, where the
$\chi_i$ range over the unramified characters.

\begin{defi}[\cite{BK1}, 4.2]
The pair $(\J,\tau)$ is said to be an $\ss$-{\it type} of $\M$ if the
irreducible representations of $\M$ whose restriction to $\J$ contains
$\tau$ are exactly those whose cuspidal support belongs to $\ss$.
\end{defi}

Thus, given an $\ss$-type $(\J,\tau)$, the functor ${\bf M}_{\tau}$
induces a bijection between the classes of irreducible representations 
of $\M$ with cuspidal support in $\ss$ and the classes of irreducible
right $\Hh(\M,\tau)$-modules.

\subsection{}
\label{RappelCovers}

Let $(\J_\M,\tau_\M)$ be an $\ss$-type of $\M$, 
and let $(\J,\tau)$ be a $\G$-cover of $(\J_\M,\tau_\M)$. 
We do not give here the definition of a cover (see \cite[8.1]{BK1}), 
which is quite technical. 
We just mention that we have $\J\cap\M=\J_\M$ and that the restriction 
of $\tau$ to $\M$ is $\tau_\M$. 
The importance of the notion of cover lies in the isomorphism 
(\ref{AApot}) below. 

Given a parabolic subgroup $\P$ of $\G$ with Levi subgroup $\M$, 
we denote by:
\begin{equation}
\label{tepe}
t_\P:\Hh(\M,\tau_\M)\to\Hh(\G,\tau)
\end{equation}
the $\CC$-alg\-ebra homomorphism given by \cite[Corollary 7.12]{BK1}.
If we denote by $\Hh$ and $\Hh_\M$ the Hecke algebras $\Hh(\G,\tau)$ 
and $\Hh(\M,\tau_\M)$, then the map $t_\P$ makes $\Hh$ into an 
$\Hh_\M$-algebra.
According to \cite{BK1} (see Theorem 8.3 and Corollary 8.4), the pair
$(\J,\tau)$ is an $\ss$-type  of $\G$ and, for any irreducible
representation $\s$ of $\M$ with cuspidal support in $\ss$, we have a
canonical $\Hh$-module isomorphism:
\begin{equation}
\label{AApot}
{\bf M}_{\tau}(\Ind_{\P}^{\G}(\s))
\simeq
\Hom_{\Hh_{\M}}(\Hh,{\bf M}_{\tau_{\M}}(\s)),
\end{equation}
where $\Ind_{\P}^{\G}$ denotes the (normalized) parabolic induction
functor.

\subsection{}
\label{GreatExpectations}
\def\UU{\U}
\def\u{u}
\def\KK{\K}
\def\KT{\tilde\K}
\def\II{\Ii}

In this paragraph, we discuss the question of the existence of types
relative to a given inertial class.
Let $\rho\in\Cc$ be a cuspidal irreducible representation, and set
$\k=\deg(\rho)$. 
Let $m$ be a positive integer which is a multiple of $k$ and let $r$
denote the positive integer such that $m=kr$.
We denote by $\ss$ the set of all multisets of size $r$ on
$\ell(\rho)$.
This is the inertial class in $\G=\G_m$ of the mult\-iset
$(\rho,\ldots,\rho)\in\M(\Cc)$, where $\rho$ occurs $r$ times. 
We have the following result:

\begin{theo}
There exists an $\ss$-type of $\G$.
\end{theo}

This is \cite[Theorem 5.5]{GSZ} if $\rho$ is of level zero (that is,
if $\rho$ has a non-zero vector invariant under the subgroup
$1+\M_k(\p_\D)$, where $\p_\D$ denotes the maximal ideal of the ring
of integers of $\D$) and \cite[Th\'eor\`eme 5.23]{VS4} if not.

\subsection{}

In order to prove Conjecture U0, we need $\ss$-types of $\G$ whose
Hecke algebras we understand precisely.
This requires the notion of simple type, which first appears in
\cite{BK} and has been generalized in \cite{VS1,VS2,VS3}.
For a def\-inition of simple type, see \cite[\S4.1]{VS3}.

\begin{prop}
\label{Estheride}
\begin{enumerate}
\item[(i)]
There is a simple type of $\G_{k}$ contained in $\rho$.
\item[(ii)]
Let $(\UU,\u)$ be a simple type contained in $\rho$.
There is a finite extension $\KK$ of $\F$ contained in
$\M_k(\D)$ such that the normalizer of $\u$ in $\G_k$ is
$\mult\KK\UU$.
\end{enumerate}
\end{prop}

\begin{proof}
Note that a type of $\G_{k}$ is contained in $\rho$ if and only if it
is a type relative to the inertial class $\ell(\rho)$.
Part (i) of the result comes from \cite[Theorem 5.4]{GSZ} if $\rho$ is
of level zero and from \cite[Th\'eor\`eme 5.21]{VS4} if not.

In order to prove part (ii), recall that the simple type $(\UU,\u)$
comes with a finite extension $\E$ of $\F$ contained in $\M_{k}(\D)$
(see \cite[\S4.1]{VS3}).
The centralizer of $\E$ in $\M_{k}(\D)$ is a central simple
$\E$-algebra isomorphic to $\M_{k'}(\D')$, where $k'$ is a positive
integer and $\D'$ a finite-dimensional central division algebra over
$\E$.
According to \cite[\S5.1]{VS3} the normalizer of $\u$ in $\G_k$ is
generated by $\UU$ and an element $\varpi$ which is a positive power
of a uniformizer of $\D'$.
The $\E$-algebra $\KK=\E[\varpi]$ is a totally ramified extension of $\E$.
As an extension of $\F$, it has the required property. 
\end{proof}

\subsection{}

In \cite[\S5.2]{VS3} one describes a process:
\begin{equation}
\label{JiPe}
(\UU,\u)\mapsto(\J,\tau)
\end{equation}
which associates, to any simple type $(\UU,\u)$ of $\G_{k}$ contained
in $\rho$, an $\ss$-type $(\J,\tau)$ of $\G$ with the following
property. 

\begin{prop}
\label{Bixiou}
For any Levi subgroup $\M$ of $\G$ containing: 
\begin{equation}
\label{DeMarsay}
\M_{0}=\G_{k}^{r}=\G_{k}\times\ldots\times\G_{k},
\end{equation}
the restriction of $(\J,\tau)$ to $\M$ is an $\ss$-type of $\M$ of
which $(\J,\tau)$ is a $\G$-cover.
\end{prop}

\begin{proof}
According to Pro\-po\-si\-tion \cite[5.5]{VS3}, the pair $(\J,\tau)$
associated to $(\UU,\u)$ by (\ref{JiPe}) is an $\ss$-type of $\G$
constructed as a cover of the type $(\UU^{r},\u^{\otimes r})$ of the
Levi subgroup $\M_{0}$.
The result follows from \cite[Proposition 8.5]{BK1}.
\end{proof}

\begin{rema}
\label{Kabbalah}
The reader should pay attention to the fact that, in general, the pair
$(\J,\tau)$ is {\it not} what we call a simple type in \cite{VS3}, but
is the type which we denote by $(\J_\P,\l_\P)$ in \cite[\S5.2]{VS3}.
Nevertheless, according to \cite[Proposition 5.4]{VS3}, there exists a
compact open subgroup $\J^{\dag}$ of $\G$ containing $\J$ such that
the induced representation of $\tau$ from $\J$ to $\J^{\dag}$ is a
simple type.
\end{rema}

\begin{exem}
\label{ExCanon}
Assume that $\D=\F$ and that $\rho$ is the trivial character of
$\GL_1(\F)$.
Then the trivial character $1_{\Oo_{\F}^{\times}}$ of the unit group
of the ring of integers $\Oo_\F$ is a simple type of $\GL_1(\F)$
containing $\rho$.
The pair $(\J,\tau)$ associated to it by (\ref{JiPe}) is the trivial
character of the standard Iwahori sub\-group of $\G=\GL_{r}(\F)$.
(By {\it standard} we mean that the reduction of $\J$ modulo $\p_\D$
is made of {\it upper} triangular matrices.) 
\end{exem}

\subsection{}
\label{MartinHeidegger}

Let $(\UU,\u)$ be a simple type contained in $\rho$ and let
$(\J,\tau)$ be the $\ss$-type of $\G$ corresponding to it by
(\ref{JiPe}).
In this paragraph, we describe the support of the Hecke algebra
$\Hh(\G,\tau)$.
Let $\KK/\F$ be as in Proposition \ref{Estheride}, let $\varpi$
be a uniformizer of $\KK$, let $\Nn$ be the normalizer of the diagonal
torus of $\GL_r(\KK)$ and let $\W$ be the subgroup of $\Nn$ made of
elements whose non-zero entries are of the form $\varpi^n$ with
$n\in\ZZ$.
As $\K$ is contained in $\Mat_k(\D)$, the group $\GL_r(\KK)$ can 
naturally be considered as a subgroup of $\G$.
Set:
\begin{equation*}
h=
\begin{pmatrix}
0&{\rm Id}_{r-1}\\
\varpi&0\\
\end{pmatrix}\in\W\subset\G,
\end{equation*}
where ${\rm Id}_{r-1}$ denotes the identity matrix of $\GL_{r-1}(\KK)$.
Note that $h$ does not normalize $\J$ in general.
According to Propositions \cite[4.3]{VS3} and \cite[5.10]{VS4}, any
element of $\Hh(\G,\tau)$ vanishes outside $\J\W\J$.
More precisely, we have the following result.

\begin{prop}
\label{RussianLessons}
Let us fix $w\in\W$.
\begin{enumerate}
\item[(i)]
The subspace of $\Hh(\G,\tau)$ made of functions supported on $\J w\J$
has dimension $1$, and any non-zero element of this subspace is
invertible.
\item[(ii)]
Let $\h\in\Hh(\G,\tau)$ be a non-zero element supported on $\J h\J$.
Then for any non-zero element $f$ supported on $\J w\J$, the
convolution product $f*\h$ (resp. $\h*f$) is supported on $\J wh\J$
(resp. on $\J hw\J$). 
\end{enumerate}
\end{prop}

\begin{proof}
We denote by $(\J^{\dag},\tau^{\dag})$ the simple type induced by
$(\J,\tau)$ (see Remark \ref{Kabbalah}).
According to \cite{VS3} (see Propositions 4.3 and 4.16 and Lemma 4.13), 
the result is true if we replace $\Hh(\G,\tau)$ by the Hecke algebra 
$\Hh(\G,\tau^\dag)$.
The result for $\Hh(\G,\tau)$ follows from 
\cite[Proposition 4.1.3 and Corollary 4.1.5]{BK}. 
\end{proof}

\begin{exem}
\label{ExCanon2}
Assume, as in Example \ref{ExCanon}, that $\D=\F$ and that $\rho$ is
the trivial character of $\GL_1(\F)$.
Then $\KK=\F$ satisfies the conditions of Proposition
\ref{Estheride}.
The choice of a uniformizer of $\F$ defines a subgroup $\W$ of
$\G=\GL_{r}(\F)$, and the Hecke algebra $\Hh(\G,\tau)$ of the trivial
character of the standard Iwahori sub\-group $\J$ of $\G$ is supported
on $\J\W\J=\G$ (the Bruhat decomposition).
\end{exem}

\subsection{}
\label{MartinLutherKing}

In this paragraph, we investigate the structure of the Hecke algebra
$\Hh(\G,\tau)$. 
Let $\KT$ be a finite {\it unramified} extension of $\KK$.
According to Examples \ref{ExCanon} and \ref{ExCanon2}, 
the trivial character $1_{\Oo_{\KT}^{\times}}$ of the unit group
of the ring of integers $\Oo_{\KT}$ is a simple type of $\GL_1(\KT)$
containing the trivial character of $\GL_1(\KT)$.
The pair associated to it by (\ref{JiPe}), which we denote by
$(\II,1_\II)$, is the trivial character of the standard Iwahori
sub\-group of $\GL_{r}(\KT)$.
Note that $\W$ can be considered as a subgroup of both $\G$ and
$\GL_{r}(\KT)$.
Given $f\in\Hh(\G,\tau)$ (resp. $f\in\Hh(\GL_r(\KT),1_\II)$), we set:
\begin{equation*}
\label{Nasser}
\supp(f)=\{w\in\W\ |\ f(w)\neq0\},
\end{equation*}
which is the support of $f$ {\it in $\W$}.
For technical reasons, this is more convenient than the support in
$\G$ (resp. in $\GL_r(\KT)$).

\begin{prop}
\label{ANenPlusFinir}
For a unique (up to isomorphism) choice of finite un\-ramified
extension $\KT$ of $\K$, there is a $\CC$-algebra isomorphism:
\begin{equation}
\label{Baudrillard*}
\Psi:\Hh(\GL_{r}(\KT),1_\II)\to\Hh(\G,\tau)
\end{equation}
such that for any function $f\in\Hh(\GL_{r}(\KT),1_\II)$, we have:
\begin{equation}
\label{SupportPreservation*}
\supp(\Psi f)=\J\cdot\supp(f)\cdot\J.
\end{equation}
\end{prop}

\begin{proof}
Theorem \cite[4.6]{VS3} gives us the result for the Hecke algebra
$\Hh(\G,\tau^\dag)$. 
The result for $\Hh(\G,\tau)$ follows from 
\cite[Proposition 4.1.3]{BK}.
\end{proof}

\begin{rema}
\begin{enumerate}
\item[(i)]
Note that (\ref{SupportPreservation*}) makes sense because $\W$ can be
seen as a subgroup of $\GL_{r}(\KT)$ on the left hand side, and of
$\G$ on the right hand side.
\item[(ii)]
The unramified extension $\KT/\K$ does not depend on the integer $r$,
but only on the cuspidal representation $\rho$.
\end{enumerate}
\end{rema}

\subsection{}
\label{MartinFiero}

Let us fix an extension $\KT$ of $\F$ as in Proposition
\ref{ANenPlusFinir}. 
Let $\P$ be the parabolic subgroup of $\G$ of upper triangular
matrices with respect to the Levi subgroup $\M_0=\G_{k}^{r}$ (see 
(\ref{DeMarsay})) and let $t_{\P}$ be the $\CC$-algebra homomorphism:
\begin{equation*}
t_\P:\Hh(\G_k^r,u^{\otimes r})\to\Hh(\G,\tau)
\end{equation*}
corresponding to $\P$ (see (\ref{tepe})).
We denote by $\Q$ the (minimal) parabolic subgroup of $\GL_r(\KT)$ of
upper triangular matrices.
Let $t_{\Q}$ be the $\CC$-algebra homomorphism:
\begin{equation*}
t_\Q:\Hh(\KT^{\times r},1_{\Oo^{\times}_{\KT}}^{\otimes r})
\to\Hh(\GL_r(\KT),1_{\II})
\end{equation*}
corresponding to $\Q$.
Let us choose a $\CC$-algebra isomorphism:
\begin{equation}
\label{Baudrillard}
\Psi_{u}:\Hh(\mult\KT,1_{\Oo^{\times}_{\KT}})\to\Hh(\G_{k},u)
\end{equation}
such that, for any function
$f\in\Hh(\mult\KT,1_{\Oo^{\times}_{\KT}})$, we have:
\begin{equation}
\label{SupportPreservation}
\supp(\Psi_u(f))=\UU\cdot\supp(f)\cdot\UU,
\end{equation}
where $\supp$ denotes the support in the group $\langle\varpi\rangle$
generated by $\varpi$, considered as a subgroup of $\mult\KT$ on the
left hand side and of $\G_{k}$ on the right hand side. 
Then there is a unique $\W$-equivariant $\CC$-algebra isomorphism:
\begin{equation*}
\label{Baudrillardr}
\Psi_{u}^{r}:\Hh(\KT^{\times r},1_{\Oo^{\times}_{\KT}}^{\otimes r})
\to\Hh(\G_{k}^{r},u^{\otimes r})
\end{equation*}
which agrees with $\Psi_{u}$ on the first tensor factor and such that, 
for any function 
$f\in\Hh(\KT^{\times r},1_{\Oo^{\times}_{\KT}}^{\otimes r})$, we have: 
\begin{equation*}
\label{SupportPreservationr}
\supp(\Psi_u^r(f))=\UU^r\cdot\supp(f)\cdot\UU^r,
\end{equation*}
where $\supp$ denotes the support in the group
$\langle\varpi\rangle^{r}$, considered as a subgroup of 
$\KT^{\times r}$ on the left hand side and of $\G_{k}^{r}$ on the
right hand side (compare \cite[7.6.19]{BK}). 
We are now ready to state the main result of this section.

\begin{theo}
\label{ExistenceHeckeMBK}
Given a $\CC$-algebra isomorphism $\Psi_u$ as in (\ref{Baudrillard}), 
there is a unique $\CC$-alg\-ebra isomorphism:
\begin{equation*}
\Psi_\G:\Hh(\GL_{r}(\KT),1_{\II})\to\Hh(\G,\tau)
\end{equation*}
such that the diagram:
\begin{equation*}
\diagram
\Hh(\GL_{r}(\KT),1_{\II})\rto^{\Psi_{\G}}&\Hh(\G,\tau)\\
\Hh(\KT^{\times r},1_{\Oo^{\times}_{\KT}}^{\otimes r})
\uto^{t_{\Q}}\rto_{\Psi_{u}^{r}}
&\Hh(\G_{k}^{r},u^{\otimes r})\uto_{t_{\P}}
\enddiagram
\end{equation*}
commutes.
\end{theo}

\begin{proof}
The proof goes {\it mutatis mutandis} as in 
\cite[Theorem 7.6.20]{BK}. 
\end{proof}

\begin{rema}
\label{Apocoloquintose}
The isomorphism $\Psi_{\G}$ preserves the canonical structure of
$\CC$-alg\-ebra with involution on the Hecke algebras (see
\S\ref{StructureUnitaire}). 
In other words, for any $f\in\Hh(\GL_{r}(\KT),1_{\II})$, we have
$\Psi_\G(f^*)=\Psi_\G(f)^*$.
This implies that unitary modules over $\Hh(\GL_{r}(\KT),1_{\II})$
correspond bijectively to unitary modules over $\Hh(\G,\tau)$.
\end{rema}


\section{Proof of Conjecture U0}

\subsection{}
\label{OnDenoting}

In this paragraph, we reduce the proof of Conjecture U0 to the
following special case:

\medskip

\noindent{\bf (S0)} \ 
{\it Let $\s,\tau\in\Irr^u$ be aligned unitary irreducible
  representations. 
Then $\s\times\tau$ is irreducible.}

\medskip

\begin{prop}
\label{AimpliesU0}
Assume that {\rm S0} holds.
Then {\rm U0} is true.
\end{prop}

\begin{proof}
Let $\s,\tau\in\Irr^u$ be irreducible unitary representations, and let:
\begin{equation*}
\s=\s_1\times\ldots\times\s_{k}
\quad\text{and}\quad
\tau=\tau_1\times\ldots\times\tau_{k'}
\end{equation*}
be the factorizations of $\s$ and $\tau$ given by Proposition
\ref{Lapounet}.
In particular, each $\s_i,\tau_j$ is simple for $1\<i\<k$ and
$1\<j\<k'$.
Moreover, we can choose the ordering such that there exists
a non-negative integer $r$ for which $\s_i$ and $\tau_i$ are 
aligned if $1\<i\<r$, and $\s_i$ is not aligned with $\tau_j$ 
if $i,j\>r+1$.
As $\s,\tau$ are unitary and irreducible, and according to
Proposition \ref{Lapounet}, each representation $\s_i,\tau_j$ is
unitary.
We write:
\begin{equation}
\label{ProJK}
\begin{split}
\s\times\tau=
(\s_1\times&\tau_1)\times\ldots\\
&\ldots\times(\s_r\times\tau_r)
\times\s_{r+1}\times\ldots\times\s_{k}\times\tau_{r+1}
\times\ldots\times\tau_{k'}.
\end{split}
\end{equation}
Assuming that S0 holds, each $\s_i\times\tau_i$ is irreducible
for $1\<i\<r$.
Therefore (\ref{ProJK}) shows that $\s\times\tau$ is a product of 
irreducible factors, no two of them being aligned.
The result now follows from Proposition \ref{Tad}.
\end{proof}

\begin{rema}
Statement S0 can be rephrased as follows: 
any simple unitary irreducible representation of a Levi 
subgroup of $\G_m$, with $m\>1$, induces irreducibly to $\G_m$. 
\end{rema}

\subsection{}
\label{U0VersionHecke}

Let $\rho\in\Cc$ be a cuspidal irreducible representation, and set
$\k=\deg(\rho)$. 
Let $m$ be a positive integer which is a multiple of $k$ and let $r$
denote the positive integer such that $m=kr$.
Let $\M$ be a Levi subgroup of $\G=\G_m$ of the form:
\begin{equation}
\label{SatLeviM}
\M=\G_{kr_1}\times\G_{kr_2},
\end{equation}
where $r_1,r_2\>1$ are positive integers such that $r_1+r_2=r$. 
Let $(\UU,u)$ be a simple type contained in $\rho$, let $(\J,\tau)$
be the $\ss$-type of $\G$ corresponding to it by (\ref{JiPe}) and let
$(\J_\M,\tau_\M)$ be the $\ss$-type of $\M$ of which $(\J,\tau)$ is a
$\G$-cover by Proposition \ref{Bixiou}. 
Let $\Hh$ and $\Hh_\M$ denote the Hecke algebras $\Hh(\G,\tau)$ and
$\Hh(\M,\tau_\M)$.
Let $\P$ be the parabolic subgroup of $\G$ of upper triangular
matrices with respect to $\M$ and let $t_{\P}$ be the $\CC$-algebra
homomorphism from $\Hh_\M$ to $\Hh$ corresponding to $\P$ 
(see (\ref{tepe})).

\begin{prop}
\label{Diderot}
Let $\V$ be a unitary irreducible $\Hh_{\M}$-module.
Then the $\Hh$-module $\Hom_{\Hh_{\M}}(\Hh,\V)$ is irreducible.
\end{prop}

\begin{proof}
We will first prove Proposition \ref{Diderot} in a particular case.
\begin{enumerate}
\item
We temporarily suppose that $\D=\F$ and that $\rho$ is the
trivial character of $\GL_1(\F)$ (see Example \ref{ExCanon}).
In that case, we can choose for $\J$ the standard Iwahori sub\-group
of $\G$ and for $\tau$ the trivial character of $\J$.
Therefore, $\J_\M$ is the standard Iwahori subgroup of $\M$ and
$\tau_\M$ is its trivial character.
The functor ${\bf M}_{\tau}$ (resp. ${\bf M}_{\tau_\M}$) associates to
a representation of $\G$ (resp. $\M$) the space of its
$\J$-invariant (resp. $\J\cap\M$-invariant) vectors.

We now recall the following crucial result of Barbasch and Moy 
\cite{BM1,BM2}.

\begin{theo}[Barbasch-Moy]
\label{JaunePenible}
The functor ${\bf M}_{\tau_\M}$ induces a bijective correspondence
between unitary irreducible representation of $\M$ with 
a non-zero space of $\J\cap\M$-invariant vectors 
and unitary irreducible right $\Hh_{\M}$-modules.
\end{theo}

Let $\s$ be an irreducible representation of $\M$ with a non-zero
space of $\J\cap\M$-invariant vectors such that 
${\bf M}_{\tau_\M}(\s)$ is isomorphic to $\V$.
By Theorem \ref{JaunePenible}, this representation is unitary. 
According to (\ref{AApot}), it is enough to prove that the
$\Hh$-module:
\begin{equation*}
{\bf M}_{\tau}(\Ind_{\P}^{\G}(\s))
=\Ind_{\P}^{\G}(\s)^{\J}
\end{equation*}
is irreducible.
According to Theorem \ref{bouloulou}, the induced representation
$\Ind_{\P}^{\G}(\s)$ is irreducible.
Because ${\bf M}_{\tau}$ preserves irreducibility, we are done.
\item
Now the symbols $\D,\rho,\J,\tau$... recover their general meaning. 
We are going to reduce the general case to our particular case 1.
Let $\KT$ be a finite extension of $\F$ as in Proposition
\ref{ANenPlusFinir}.
We use the notations of
\S\S\ref{MartinLutherKing}--\ref{MartinFiero}.
Let $\L$ denote the Levi sub\-group:
\begin{equation*}
\L=\GL_{r_1}(\KT)\times\GL_{r_2}(\KT).
\end{equation*}
Let $\Q$ be the parabolic subgroup of $\GL_r(\KT)$ of upper triangular
matrices with respect to $\L$ and let $t_{\Q}$ be the $\CC$-algebra
homomorphism from the Hecke algebra
$\Hh_\L=\Hh(\L,1_{\II\cap\L})$ to $\Hh(\GL_r(\KT),1_{\II})$
corresponding to $\Q$.
Let $\Psi_\G$ denote the $\CC$-algebra isomorphism of Theorem
\ref{ExistenceHeckeMBK}. 

\begin{prop}
\label{ExistenceHecke}
There is a $\CC$-algebra isomorphism:
\begin{equation*}
\Psi_\M:\Hh(\L,1_{\II\cap\L})\to\Hh(\M,\tau_\M)
\end{equation*}
such that the diagram:
\begin{equation*}
\diagram
\Hh(\GL_{r}(\KT),1_{\II})\rto^{\Psi_{\G}}&\Hh(\G,\tau)\\
\Hh(\L,1_{\II\cap\L})\uto^{t_{\Q}}\rto_{\Psi_\M}
&\Hh(\M,\tau_\M)\uto_{t_{\P}}
\enddiagram
\end{equation*}
commutes.
\end{prop}

\begin{proof}
According to Theorem \ref{ExistenceHeckeMBK}, it suffices to choose
for $\Psi_\M$ the $\W$-equivariant $\CC$-algebra isomorphism 
which agrees with $\Psi_{\G_{kr_1}}$ on the first tensor factor and
such that we have:
\begin{equation*}
\supp(\Psi_\M(f))=\J_\M\cdot\supp(f)\cdot\J_\M
\end{equation*}
for any function $f\in\Hh(\L,1_{\II\cap\L})$.
\end{proof}

This allows us to make $\V$ into a module over $\Hh_\L$, and thus to
identify the $\Hh$-module $\Hom_{\Hh_{\M}}(\Hh,\V)$ with the 
$\Hh(\GL_{r}(\KT),1_{\II})$-module given by:
\begin{equation}
\label{Ekatarina}
\Hom_{\Hh_\L}(\Hh(\GL_r(\KT),1_{\II}),\V).
\end{equation}
As $\Psi_\M$ preserves the canonical structure of $\CC$-alg\-ebra with
involution (see Remark \ref{Apocoloquintose}), $\V$ is irreducible and
unitary as a $\Hh_\L$-module.
Therefore (\ref{Ekatarina}) is irreducible according to case 1.
\end{enumerate}

This ends the proof of Proposition \ref{Diderot}.
\end{proof}

\subsection{}

In this paragraph, we prove S0.
With the notations of \S\ref{U0VersionHecke}, it suffices to prove the 
following result. 

\begin{prop}
\label{BAlaBA}
Let $\s$ be a simple unitary irreducible representation of $\M$ with
cuspidal support in $\M(\ell(\rho))$.
Then the induced representation $\Ind_\P^{\G}(\s)$ is irreducible.
\end{prop}

\begin{proof}
We apply Proposition \ref{Diderot} to the 
$\Hh_{\M}$-module $\V={\bf M}_{\tau_\M}(\s)$, which is irreducible and
unitary (see \S\ref{StructureUnitaire}).
The $\Hh$-module ${\bf M}_{\tau}(\Ind_\P^{\G}(\s))$ is then
irreducible, thanks to (\ref{AApot}).
The result now follows from the fact that ${\bf M}_{\tau}$ preserves
reducibility.
\end{proof}

This ends the proof of Conjecture U0, thanks to Proposition
\ref{AimpliesU0}. 

\begin{rema}
\label{ArbiChar}
In \cite{Tadic2}, as in this paper, the characteristic of $\F$ is
assumed to be zero.
However, with the works of Badulescu \cite{Bathese,Ba2} and 
M\'\i nguez \cite{Minguez}, this assumption seems to be superfluous,
and the Tadi\'c class\-ification of the unitary dual of $\GL_m(\D)$
should be available in arbitrary characteristic.
More precisely, when $\F$ is of positive characteristic:
\begin{enumerate}
\item[(1)]
M\'\i nguez \cite[\S2.1.14]{Minguez} proved that the ring $\R$ of
\S\ref{UZero} is commutative;
\item[(2)] 
Badulescu \cite{Ba2} proved that any square integrable 
irreducible representation of a Levi subgroup of $\G_m$ induces
irreducibly to $\G_m$ (see \S\ref{JJR}).
\end{enumerate}
It would therefore be interesting to write down a classification 
of the unitary dual of $\GL_m(\D)$ with no assumption on the
characteristic of $\F$.
\end{rema}


\section{Reducibility points}
\label{Sec4}

Let $\rho\in\Cc$ be a cuspidal irreducible representation of degree
$k$. 
In this section, we determine the unramified characters $\chi$ of
$\G_k$ such that the representation $\rho\times\rho\chi$ is reducible. 
This could provide a definition of the integer $b(\rho)$ of
\S\ref{SerieDiscrete} without refering to the Jacquet-Langlands
correspondence.

\subsection{}

Let $(\U,\u)$ be a simple type contained in $\rho$.
According to Proposition \ref{Estheride}, the norm\-alizer $\N$ of
$\u$ in $\G_\k$ is generated by $\U$ and a uniformizer $\varpi$ of 
the extension $\K$.
Let $q_\F$ denote the cardinal of the residue field of $\F$.

\begin{prop}
\label{Vandale}
The group of unramified characters $\chi$ of $\G_k$ such that
$\rho\simeq\rho\chi$ is finite. 
\end{prop}

\begin{proof}
According to \cite[\S5.1]{VS3}, the representation $\u$ extends to
an irreducible representation $\widetilde\u$ of $\N$, and $\rho$ is
equivalent to the representation of $\G_k$ compactly induced from
$\widetilde\u$.
Given an unramified character $\chi$ of $\G_k$, the representation
$\rho\chi$ is compactly induced from the restriction 
$\widetilde\u\chi_{|\N}$, which
is equivalent to $\widetilde\u$ if and only if $\chi(\varpi)=1$.
Let us define a positive integer $n$ by:
\begin{equation}
\label{Epanadiplose}
\nu(\varpi)=q_\F^{-n}.
\end{equation}
Then the group of unramified characters $\chi$ of $\G_k$ such that
$\rho\simeq\rho\chi$ is cyclic of order $n$.
\end{proof}

\begin{defi}
The {\it torsion number} of $\rho$, which we denote by $n(\rho)$, 
is the cardinal of the group of unramified characters $\chi$ of 
$\G_k$ such that $\rho\simeq\rho\chi$.
\end{defi}

\subsection{}

Let $\h$ be a non-trivial element of the Hecke algebra $\Hh(\G_k,\u)$
supported by the double coset $\U\varpi\U$ (which actually is a single
coset).
According to Propositions \ref{RussianLessons} and
\ref{ANenPlusFinir}, such an element is invertible and $\Hh(\G_k,\u)$
is the commutative $\CC$-algebra generated by $\h$ and $\h^{-1}$.
Therefore, the irreducible $\Hh(\G_k,\u)$-modules are one-dimensional
and characterised, up to isomorphism, by a non-zero complex number
given by the eigenvalue of $\h$.

\begin{defi}
If $\V$ is an irreducible $\Hh(\G_k,\u)$-module on which $\h$ acts by
$\l\in\mult\CC$ and $\chi$ an unramified character of $\G_k$, we will
denote by $\V\chi$ the irreducible $\Hh(\G_k,\u)$-module (with the
same underlying space as $\V$) on which $\h$ acts by $\chi(\varpi)\l$. 
\end{defi}

Let ${\bf M}={\bf M}_{\u}$ denote the functor defined by
(\ref{FoncMT}) relative to the pair $(\U,\u)$.
It induces a bijective correspondence between the inertial class
$\ell(\rho)$ of $\rho$ and the set of all classes of irreducible
$\Hh(\G_k,\u)$-modules.

\begin{lemm}
\label{Fulgence}
For any unramified character $\chi$ of $\G_k$, the module 
${\bf M}(\rho\chi)$ is equal to ${\bf M}(\rho)\chi^{-1}$.
\end{lemm}

\begin{proof}
This is proved in \cite[\S2]{BK1}.
The reader should pay attention to the fact that in \cite{BK1}, 
the symbol $\Hh(\G_k,\u)$ has a slightly different meaning.
To recover our $\Hh(\G_k,\u)$, one has to apply the isomorphism 
given by \cite[(2.3)]{BK1}.
\end{proof}

Let $(\J,\tau)$ be the type of $\G_{2k}$ which corresponds to
$(\U,\u)$ by (\ref{JiPe}).
This is a $\G_{2k}$-cover of the pair $(\U^2,\u^{\otimes2})$
considered as a type of the Levi sub\-group $\M=\G_k\times\G_k$, so
that we have $(\J_\M,\tau_\M)=(\U^2,\u^{\otimes2})$.
Let $\Hh$ and $\Hh_\M$ denote the Hecke algebras relative to $\tau$
and $\tau_\M$ respectively.
Let ${\bf M}_{\tau}$ be the functor which corresponds to $\tau$, let
$\P$ be the parabolic subgroup of $\G_{2k}$ of upper triangular
matrices relative to $\M$ and let $t_\P$ be the map given by
(\ref{tepe}).
Let $\KT$ be a finite extension of $\F$ as in Proposition
\ref{ExistenceHecke} and let $q_{\KT}$ be the cardinal of its residue
field.

\begin{prop}
\label{Ashk}
Let $\V$ be an irreducible $\Hh(\G_k,\u)$-module and let $\chi$ be an
unramified character of $\G_k$.
Then the $\Hh$-module:
\begin{equation}
\label{Sepharade}
\Hom_{\Hh_{\M}}(\Hh,\V\otimes\V\chi^{-1})
\end{equation}
is reducible if and only if $\chi(\varpi)=q_{\KT}$ or
$\chi(\varpi)=q_{\KT}^{-1}$. 
\end{prop}

\begin{proof}
Let $\s$ be the unramified twist of $\rho$ such that ${\bf M}(\s)$ is
isomorphic to $\V$.
According to Lemma \ref{Fulgence}, we have a canonical $\Hh$-module
isomorphism: 
\begin{equation}
\label{EulalalilalaProto}
{\bf M}_{\tau}(\s\times\s\chi)\simeq
\Hom_{\Hh_\M}(\Hh,\V\otimes\V\chi^{-1}).
\end{equation}
\begin{enumerate}
\item
We temporarily suppose that $\D=\F$ and that $\rho$ is the trivial
character of $\GL_1(\F)$.
In that case, we can choose for $\U$ the maximal compact subgroup of
$\mult\F$ and for $\u$ the trivial character of $\U$.
We have $n(\rho)=1$ and $\KT=\F$, and the representation
$\s\times\s\chi$ is reducible if and only if $\chi=|\ |_\F$ or
$\chi=|\ |_\F^{-1}$.
\item
Let $\Ii$ denote the standard Iwahori subgroup of $\GL_2(\KT)$ and
$1_\Ii$ its trivial character, which is the $\GL_2(\KT)$-cover
associated by (\ref{JiPe}) to the trivial character, which we 
denote by $1_{\mult\Oo_{\KT}}$, of the maximal compact subgroup of
$\mult\KT$.
Let $\L$ denote the Levi sub\-group $\GL_1(\KT)\times\GL_1(\KT)$,
let $\Q$ be the parabolic subgroup of $\GL_2(\K)$ of upper triangular
matrices relative to $\L$ and let $t_{\Q}$ be the $\CC$-algebra
homomorphism from $\Hh_\L=\Hh(\L,1_{\Ii\cap\L})$ to
$\Hh(\GL_2(\KT),1_\Ii)$ corresponding to $\Q$.

We make $\V$ into a module over $\Hh(\mult\KT,1_{\mult\Oo_{\KT}})$ by
fixing a $\CC$-algebra isomorphism (\ref{Baudrillard}), 
which allows us, according to Proposition \ref{ExistenceHecke}, to
identify the $\Hh$-module (\ref{Sepharade}) with the
$\Hh(\GL_2(\KT),1_\Ii)$-module:
\begin{equation}
\label{NomDeLaRose}
\Hom_{\Hh_{\L}}(\Hh(\GL_2(\KT),1_\Ii),\V\otimes\V\tilde\chi^{-1}),
\end{equation}
where $\tilde\chi$ denotes the unramified character of $\mult\KT$
which takes the same value as $\chi$ on $\varpi$.
According to case 1, this module is reducible if and only if 
$\tilde\chi=|\ |_{\KT}$ or $\tilde\chi=|\ |_{\KT}^{-1}$, which amounts
to saying that (\ref{NomDeLaRose}) is reducible if and only if
$\chi(\varpi)=q_{\KT}$ or $\chi(\varpi)=q_{\KT}^{-1}$.
\end{enumerate}
This gives us the required result.
\end{proof}

Let $f(\rho)$ denote the residue degree of $\KT$ over $\F$.
We state the main result of this section.

\begin{theo}
\label{MainTheo4}
Let $s\in\CC$.
Then $\rho\times\nu^{s}\rho$ is reducible if and only if: 
\begin{equation*}
s=f(\rho)n(\rho)^{-1}\quad\text{or}\quad s=-f(\rho)n(\rho)^{-1}.
\end{equation*}
\end{theo}

\begin{proof} 
We apply Proposition \ref{Ashk} with the unramified character 
$\chi=\nu^s$.
The result follows from the definition of $n(\rho)$ by
(\ref{Epanadiplose}).
\end{proof}



\providecommand{\bysame}{\leavevmode ---\ }
\providecommand{\og}{``}
\providecommand{\fg}{''}
\providecommand{\smfandname}{\&}
\providecommand{\smfedsname}{\'eds.}
\providecommand{\smfedname}{\'ed.}
\providecommand{\smfmastersthesisname}{M\'emoire}
\providecommand{\smfphdthesisname}{Th\`ese}

\end{document}